\documentclass[12pt]{article}
\usepackage{amssymb,amsmath}
\usepackage{graphics}
\newtheorem{thm}{Theorem}
\newtheorem{cor}[thm]{Corollary}
\newtheorem{lemma}[thm]{Lemma}
\newenvironment{defin}{\medskip\noindent{\sc
Definition}.}{\goodbreak\medskip}
\newenvironment{nota}{\medskip\noindent{\sc
Notations}.}{\goodbreak\medskip}
\newenvironment{remk}{\noindent{\sc
Remark}.}{\goodbreak\vskip10pt}
\newtheorem{prop}[thm]{Proposition}

\def\demo{\medskip\goodbreak\noindent
     \hbox{\sc Proof \kern .3em}\ignorespaces}%
  \def \qedbox{$\square$}%
  \def \qed{\hglue1mm\hfill{\ifmmode\qedbox
     \else\unskip\ \hglue0mm\hfill\qedbox\medskip
      \goodbreak\fi}}%
\def\enddemo{\qed\goodbreak\vskip10pt}%

\def\qed{\hglue1mm\hfill\raise -2pt\hbox{\vrule\vbox to 10pt{\hrule width
4pt
                  \vfill\hrule}\vrule}}

\newcommand{\T}{\mathbb {T}}

\newcommand{\esse}{\mathbb {S}}
\newcommand{\R}{\mathbb {R}}

\newcommand{\Ic}{\mathcal {I}}

\newcommand{\Mc}{\mathcal {M}}

\newcommand{\Gc}{\mathcal {G}}
\newcommand{\Lc}{\mathcal {L}}

\newcommand{\Ac}{\mathcal {A}}
\newcommand{\Tc}{\mathcal {T}}

\begin{document}
\title{{ On a theorem due to Birkhoff}}
\author{M.-C. ARNAUD
\thanks{ANR Project BLANC07-3\_187245, Hamilton-Jacobi and Weak KAM Theory}
\thanks{ANR DynNonHyp}
\thanks{Universit\'e d'Avignon et des Pays de Vaucluse, Laboratoire d'Analyse non lin\' eaire et G\' eom\' etrie (EA 2151),  F-84 018Avignon,
France. e-mail: Marie-Claude.Arnaud@univ-avignon.fr}
}
\maketitle
\abstract{   The manifold $M$ being closed and connected, we prove that every submanifold of $T^*M$ that is Hamiltonianly isotopic to the zero-section and that is invariant by a Tonelli flow is a graph.
}
\newpage
\tableofcontents
\newpage
\section{Introduction}
A famous theorem due to G.~D.~Birkhoff asserts that any essential invariant curve that is invariant by an area preserving  twist map of the annulus is the graph of a continuous map (see \cite{Bir1}, \cite{He1}, \cite{Fa2}, \cite{K-O}, \cite{Sib}). Since that, a lot of attempts were made to generalize this result to higher dimensions. Under some assumptions, the authors prove that for a convex Hamiltonian of a cotangent bundle or a multidimensional positive twist map,  an invariant Lagrangian manifold that is Hamiltonianly isotopic to the zero section is a graph. In   general, the hypothesis is that the dynamic restricted to the invariant manifold is chain recurrent (see \cite{Bia1}, \cite{BiaPol2}, \cite{BiaPol1}, \cite{He2}, \cite{BiaPol3}  ...). In \cite{BiaPol3}, the authors ask if the result is true without such an assumption and say~: ``we have neither a proof nor a counterexample.''\\
 We will see that in the case of a Tonelli Hamiltonian, this hypothesis is useless. We will prove~:
 
 \begin{thm}\label{th1}
 Let $M$ be a compact and connected manifold. Let $H~: T^*M\rightarrow \R$ be a Tonelli Hamiltonian and let $\Tc\subset T^*M$ be an invariant Lagrangian submanifold that is Hamiltonianly isotopic to the zero section. Then $\Tc$  is in fact a Lagrangian graph.
 \end{thm}
 
 The main argument to prove this theorem is the so-called weak KAM theory. This theory was built in the 90's by A.~Fathi (see \cite{Fa1}). Another important ingredient of the proof is the use of a graph selector, or more exactly a function selector. These graph selectors give us a way to choose a pseudograph (it is a kind of discontinuous exact Lagrangian graph) in the initial Lagrangian manifold; they were firstly introduced by M.~Chaperon (see \cite{Cha}) and Y.~Oh (see \cite{Oh}); we will use the construction given by G.~Paternain, L.~Polterovich and  K.~Siburg in \cite{Patpoltsib}; in this paper, a very interesting comparison is done between the graph selector and some weak KAM solutions and we will go on with this comparison.  Let us mention too the preprint \cite{O-V} of A.~Ottolenghi and  C. Viterbo where a construction is given of the so-called ``geometric solution of the Hamilton-Jacobi equation.'' In this last paper too, the authors compare their geometric solution with other solutions, the viscosity ones (that are equal to the weak KAM ones in the autonomous case), but their result is valuable only for the time-dependent case. Curiously, we will prove that in the autonomous case too, the geometric solution corresponds to a weak KAM solution $u$. As $u$ is  a positive and negative weak KAM solution, we will deduce that the initial Lagrangian submanifold is in fact the graph of $du$. 
 \section{ A function selector}\label{select}
Let us recall the construction of a graph selector that is given in  \cite{Patpoltsib}. If $N\subset T^*M$ is a Lagrangian submanifold that is Hamiltonianly isotopic to the zero section, we can associate with it (see \cite{Sik1}) what is called a 
generating 
function quadratic at infinity (gfqi) where~:

\begin{defin} If $N\subset T^*M$ is Lagrangian, it admits a gfqi $S$ if there exists a smooth function $S~: (q, \zeta)\in M\times \R^N\rightarrow S(q, \zeta)\in \R$ such that~:\\
$\bullet$  $0$ is a regular value of the map $\frac{\partial S}{\partial\zeta}$; we introduce the notation~: $\Sigma_S=\{ (q, \zeta)\in M\times \R^N; \frac{\partial S}{\partial \zeta}(q, \zeta)=0\}$; then $\Sigma_S$ is a submanifold of $M\times \R^N$ that has the same dimension as $M$;\\
$\bullet$ a compact set $K\subset M\times \R^N$ exists so that, for every $q\in M$, the restriction of $S$ to $(\{q\}\times \R^N)\backslash K$ is a non-degenerate quadratic form;\\
$\bullet$ the map $i_S~: \Sigma_S\rightarrow T^*M$ defined by $i_S(q, \zeta)=(q, \frac{\partial S}{\partial q}(q, \zeta))$ is an embedding such that $i_S(\Sigma_S)=N$.

In this case, we have~: $$N=\{ (q, d_qS(q, \zeta)); d_\zeta S(q, \zeta)=0\}.$$
\end{defin}
Such a generating function is used in \cite{Patpoltsib} to construct a Lipschitz function $\Phi~: M\rightarrow \R$ via a min-max method. This Lipschitz function $\Phi$ satisfies~: \\
$\bullet$ for all $q\in M$, $\Phi(q)$ is a critical value of $S(q, .)$;\\
$\bullet$ there exists a dense open subset $U_0$ of $M$ with full Lebesgue measure such that $\Phi$ is differentiable on $U_0$ and~: $\forall q\in U_0, (q, d\Phi(q))\in N$. Moreover~: $\forall q\in U_0, \Phi(q)=S\circ i_S^{-1}(q, d\Phi (q))$.

In \cite{Patpoltsib}, the function is called a ``graph selector'', because it is used to select a part of the initial Lagrangian manifold $N$~: $\{ q, d\Phi(q)); q\in U_0\}\subset N$. But this function is  more than just a graph selector~: in fact, the function $\Phi$ is a means of selecting a value  of $S$ above every point $q\in M$. This is important because in the weak KAM formalism, we use   continuous  functions and not just  discontinuous Lagrangian graphs.

\bigskip
We will prove in section \ref{sectionproof} that if $N$ is invariant by a Tonelli flow, then $\Phi$ is a $C^1$ function. In this case,  the graph of $d\Phi$ is a submanifold of $N$ that has the same dimension as $N$. As $N$ is connected (because $M$ is), then $N$ is the graph of the $C^0$ map $d\Phi$.  A classical result asserts that is the $C^0$ graph of $d\Phi$ is invariant by a Tonelli flow, then $d\Phi$ is Lipschitz. Being a smooth manifold that is the graph of a Lipschitz function, $N$ is then the graph of the smooth function $d\Phi$.

   \section{Weak KAM theory} Except proposition \ref{propkam} and its corollary, all the results of this section are proved in \cite{Fa1} or \cite{Be1}.\\
   Let us recall that a Tonelli Hamiltonian is a $C^3$ function $H~: T^*M\rightarrow \R$ that is~:\\
$\bullet$ superlinear in the fiber~: $\forall A\in \R, \exists B\in \R, \forall (q,p)\in T^*M, \| p\|\geq B\Rightarrow  H(q,p) \geq A\| p\|$;\\
$\bullet$ $C^2$-convex in the fiber~: for every $(q,p)\in T^*M$, the Hessian $\frac{\partial^2H}{\partial p^2}$ of $H$ in the fiber direction is positive definite as a quadratic form.\\
We denote the Hamiltonian flow of $H$ by $(\varphi_t)$ and the Hamiltonian vector-field by $X_H$. \\
A Lagrangian function $L~: TM\rightarrow \R$ is associated with $H$. It is defined by~:
$\displaystyle{L(q, v)=\max_{p\in T^*_qM} (p.v-H(q,p))}$. Then $L$ is $C^2$-convex and superlinear in the fiber and has the same regularity as $H$. We denote its Euler-Lagrange flow by $(f_t)$. Then $(\varphi_t)$ and $(f_t)$ are conjugated by the Legendre map~: $\Lc~: (q, p)\in T^*M\rightarrow (q, \frac{\partial H}{\partial p}(q, p))\in TM$; more precisely, we have~: $\Lc\circ \varphi_t=f_t\circ \Lc$.
\subsection{Domination property}
\subsubsection{Semigroups of Lax-Oleinik}
Following A.~Fathi (see \cite{Fa1}), we may associate two semi-groups, called Lax-Oleinik semi-groups, to any Tonelli Hamiltonian~:\\
$\bullet$ the negative Lax-Oleinik semi-group $(T_t^-)_{t>0}$ is defined by~:
$$\forall u\in C^0(M, \R),  T_t^-u(q)=\min_{q'\in M} \left(u(\gamma(0))+\int_0^tL(\gamma(s), \dot\gamma(s))ds\right);$$
where the infimum is taken on the set of $C^1$ curves  $\gamma~: [0, t]\rightarrow M$ such that   $\gamma (t)=q $. \\
$\bullet$  the positive Lax-Oleinik semi-group is defined by~:
$$\forall u\in C^0(M, \R),  T^+_tu(q)=\max_{q'\in M} \left(u(\gamma(t))-\int_0^tL(\gamma (s), \dot\gamma(s))ds\right);$$
where the infimum is taken on the set of $C^1$ curves  $\gamma~: [0, t]\rightarrow M$ such that $\gamma (0)=q $. 
\subsubsection{Dominated functions}
If $u\in C^0(M, \R)$ and $k\in \R$, we write $u\prec L+k$ and we say that $u$
 is dominated by $L+k$ if for each $C^1$ curve $\gamma~: [a, b]\rightarrow M$, we have~:
 $$u(\gamma(b))-u(\gamma (a))\leq \int_a^bL(\gamma (s), \dot\gamma (s))ds +k(b-a).$$
 Then, we have~: $u\prec L+k\Leftrightarrow \forall t\geq 0, u\leq kt+T_t^-u\Leftrightarrow \forall t\geq 0, T_t^+u-kt\leq u$.\\
 It is proved in \cite{Fa1} that such a dominated function $u$  is Lipschitz, hence differentiable almost everywhere  and satisfies~: $H(q, du(q))\leq k$ at every point $q$ of $M$ where $u$ is differentiable. Moreover, it is proved too that every Lipschitz function $u~: M\rightarrow \R$ such that at Lebesgue almost every point $q$,   $u$ is differentiable and ~: $H(q, du(q))\leq k$, is dominated by $L+k$.
 \subsubsection{weak KAM solutions and Ma\~n\'e's critical value}
A function $u~: M\rightarrow \R$ is a negative (resp. positive) weak KAM solution if there exists $c\in\R$ such that~: $\forall t>0, T_t^-u=u-ct$ (resp. $\forall t>0,   T_t ^+u=u+ct$). Then there   exists  at least one positive and one negative   weak K.A.M. solutions  (see \cite{Fa1} or \cite{Be1}).  The constant $c$ is   unique and is called Ma\~n\'e's critical value. 
  
Many characterizations of Ma\~n\'e's critical value exist. For example, it is proved in \cite{CIPP1} that~:
$$c=\inf_{u\in C^\infty (M, \R)}\max_{q\in M}H(q, du(q)).$$
 
 Ma\~n\'e's critical value is   the greatest lower bound of the set of the numbers $k\in \R$ for which there exists $u\in C^0(M, \R)$ with $u\prec L+k$. 
 
 An interesting property of the weak WAM solutions is the forward (resp. backward) invariance of their pseudographs. If $u~: M\rightarrow \R$ is a Lipschitz function, we denote the graph  of $du$ by $\Gc (du)$~: 
 $\Gc(du)=\{ (q, du(q)); u \quad {\rm is}\quad {\rm differentiable}\quad{\rm at}\quad q\}$. Then, if $u_-$ (resp. $u_+$) is a negative (resp. positive) weak KAM solution, we have~: $\forall t>0, \varphi_t(\overline{\Gc(du_-)})\subset \Gc(du_-)$ (resp. $\varphi_{-t}(\overline{\Gc(du_+)})\subset \Gc(du_+)$).
 \subsection{Mather set, Aubry set and Peierls barrier}
 \subsubsection{Minimizing orbits and measures}
 Let us introduce a notation~: 

 \begin{nota}
  If $t>0$, the function $A_t~: M\times M\rightarrow \R$ is defined by~: 
$$A_t(q_0, q_1)=\inf_\gamma \int_0^tL(\gamma (s), \dot\gamma (s))ds=\min_\gamma \int_0^tL(\gamma (s), \dot\gamma (s))ds$$
where the infimum is taken on the set of $C^1$ curves  $\gamma~: [0, t]\rightarrow M$ such that $\gamma (0)=q_0$ and $\gamma (t)=q_1$.  Let us recall that $\gamma_0~: [0, t]\rightarrow M$ is a critical point of $A_t$ on the set of $C^1$ curves  $\gamma~: [0, t]\rightarrow M$ such that $\gamma (0)=q_0$ and $\gamma (t)=q_1$ if, and only if, $(\gamma, \dot\gamma)$ is an orbit piece for the Euler-Lagrange flow.  We say that $\gamma$ is minimizing if it achieves the minimum in the previous equality. Moreover, $\gamma~: \R\rightarrow M$ is minimizing if its restriction to every segment is minimizing. The corresponding orbits (for $(f_t)$ and $(\varphi_t)$) are said to be minimizing. An invariant Borel probability measure with compact support is said to be minimizing if its support is filled with minimizing orbits.
 \end{nota}
\subsubsection{Mather set and conjugate weak KAM solutions}
Let us introduce the Mather set~:

\begin{defin}
The Mather set, denoted by $\Mc^*(H)$,  is the union of the supports of the minimizing measures. The projected Mather set is $\Mc(H)=\pi(\Mc^*(H))$ where $\pi~: T^*M\rightarrow M$ is the projection.
\end{defin}
J.~Mather proved in \cite{mather1} that $\Mc^*(H)$ is compact, non-empty and that it is a Lipschitz graph above a compact part of the zero-section of $T^*M$.

A.~Fathi proved in \cite{Fa1} that if $u_-$ is a negative weak KAM solution, then there exists a unique positive weak KAM solution $u_+$ such that $u_{-|\Mc(H)}=u_{+|\Mc(H)}$. Such a pair $(u_-, u_+)$ is called a pair of conjugate weak KAM solutions. For such a pair, we have~:
  \begin{enumerate}
  \item[$\bullet$] $\forall q\in \Mc(H), u_-(q)=u_+(q)$; let us denote the set of equality~: $\Ic (u_-, u_+)=\{ q; u_-(q)=u_+(q)\}$  by $\Ic (u_-,u_+)$; then $\Mc(H)\subset \Ic(u_-, u_+)$; 
  \item[$\bullet$] $u_-$ and $u_+$ are differentiable at every point $q\in \Ic(u_-, u_+)$;   when $q\in\Mc (H)$ and  $(q, p)\in  \Mc^*(H)$ is its lift to $\Mc^*(H)$, then $  du_-(q)=  du_+(q)=p$;
  \item[$\bullet$] $u_+\leq u_-$.
  \end{enumerate}
  Moreover, if $u~: M\rightarrow \R$ is a function such that $u\prec L+c$, then there exists a unique pair $(u_-, u_+)$ of conjugate weak KAM solutions such that $u_{-|{\cal M}(H)}=u_{+|{\cal M}(H)}=u_{|{\cal M}(H)}$. In this case, we have~: $u_+\leq u\leq u_-$.
 \subsubsection{Aubry set}
If $(u_-, u_+)$ is a pair of conjugate weak KAM solutions, we denote by ${\cal I}(u_-, u_+)$ the set of equality~:
$${\cal I}(u_-, u_+)=\{ q\in M; u_-(q)=u_+(q)\}.$$
Then ${\cal M}(H)\subset {\cal I}(u_-, u_+)$, the two functions $u_-$ and $u_+$ are differentiable at every point of ${\cal I}(u_-, u_+)$ and their derivatives are equal on this set. We denote by $\tilde{\cal I}(u_-, u_+)$ the following lift of ${\cal I}(u_-, u_+)$~: 
$$\tilde{\cal I}(u_-, u_+)=\{ (q, du_-(q)); q\in {\cal I}(u_-, u_+)\} = \{ (q, du_+(q)); q\in {\cal I}(u_-, u_+)\} .$$
We have~: ${\cal M}^*(H)\subset\tilde{\cal I}(u_-, u_+)$ and it is proved in \cite{Fa1} that $\tilde{\cal I}(u_-, u_+)$ is a Lipschitz graph above ${\cal I}(u_-, u_+)$.\\
The Aubry set is defined by~:
$${\cal A}^*(H)=\bigcap\tilde{\cal I}(u_-,u_+)$$
where the intersection  is taken on the set of pairs $(u_-, u_+)$ of conjugate weak KAM solutions. The projected Aubry set is~: ${\cal A}(H)=\pi({\cal A}^*(H))$. Then ${\cal A}^*(H)$ is a Lipschitz graph above ${\cal A}(H)$ that is closed, non-empty and invariant.

The Peierls barrier   $h~: M\times M\rightarrow \R$ is defined by~: $\displaystyle{h(q_1, q_2)=\liminf_{T\rightarrow +\infty}(A_T(q_1, q_2)+cT)}$. It is proved in \cite{Fa1} that $h$ is Lipschitz and that the previous $\liminf$ is in fact a true limit, and even an uniform limit. Moreover, we have~:
\begin{enumerate}
\item[$\bullet$] for every $q\in M$~: $h(q,q)\geq 0$;
\item[$\bullet$] if $u\prec L+c$, then~: $\forall q_1, q_2\in M, u(q_2)-u(q_1)\leq h(q_1, q_2)$;
\item[$\bullet$] $\forall q\in M, q\in {\cal A}(H)\Leftrightarrow h(q, q)=0$.
\end{enumerate}
We deduce easily that $q$ belongs to $\Ac(H)$ if and only if there exists a sequence $(t_n)\in\R_+$ tending to $+\infty$ and a sequence of curves $\gamma_n~: [0, t_n]\rightarrow M$ such that $\displaystyle{\lim_{n\rightarrow \infty} \gamma_n(0)=\lim_{n\rightarrow \infty}\gamma_n(t_n)=q_n}$ and 
 $\displaystyle{\lim_{n\rightarrow \infty}\int_0^{t_n}(L(\gamma_n, \dot\gamma_n)+c)\leq 0}$. In this case, the last limit is equal to $0$ and\\
  $\displaystyle{\lim_{n\rightarrow \infty}(\gamma_n(0),\frac{\partial L}{\partial v}(\gamma_n(0), \dot\gamma_n(0)))=\lim_{n\rightarrow \infty}(\gamma_n(t_n),\frac{\partial L}{\partial v}(\gamma_n(t_n), \dot\gamma_n(t_n)))}$ is the point of $\Ac^*(H)$ that is above $q$.
  \subsection{More on the weak KAM theory}
  In \cite{Fa1}, Albert Fathi proves that a function that is a positive and negative weak KAM solution is $C^{1,1}$. Let us now give a result that may be useful to  prove that some functions are positive and negative weak KAM solutions.
  
\begin{prop}\label{propkam}
Let $u~: M\rightarrow \R$ be a dominated function~: $u\prec L+c$ and $(u_-, u_+)$ the pair of conjugate weak KAM solutions such that $u=u_-=u_+$. Then~:
\begin{enumerate}
\item[$\bullet$] if~: for almost $q\in M, \exists q_0\in {\cal A}(H), u(q_0)-u(q)\geq h(q , q_0)$, then $u=u_+$;
\item[$\bullet$]  if~: for almost $q\in M, \exists q_0\in {\cal A}(H), u(q)-u(q_0)\geq h(q_0, q)$, then $u=u_-$.
\end{enumerate}
\end{prop}
\demo
We only prove the first point, the second one being similar.  We know that $u_+\leq u\leq u_-$ and that~: $\forall q_0\in \Ac(H), u(q_0)=u_-(q_0)= u_+(q_0)$. \\
Let us now consider $q\in M$ such that there exists $q_0\in \Ac(H)$ such that $u(q_0)-u(q)\geq h(q, q_0)$. As $u_+\leq u$ and $u(q_0)=u_+(q_0)$, we have~: $u(q_0)-u(q)\leq u_+(q_0)-u_+(q)$. As $u_+$ is a weak KAM solution, it is dominated by $L+c$. We have then~:
$$h(q, q_0)\leq u(q_0)-u(q)\leq u_+(q_0)-u_+(q)\leq h(q, q_0).$$
We deduce that $ u(q_0)-u(q)= u_+(q_0)-u_+(q)$ and then $u(q)=u_+(q)$. The two functions $u$ and $u_+$ are   continuous and equal almost everywhere, they are then equal everywhere.
\enddemo
\begin{cor}\label{corkam}
Let $u~: M\rightarrow \R$ be a dominated function~: $u\prec L+c$   such that~:  for almost  $q\in M, \exists q_1, q_2 \in {\cal A}(H), u(q_1)-u(q)\geq h(q , q_1)$ and $u(q)-u(q_2)\geq h(q_2, q)$.
Then $u$ is $C^{1, 1}$ and the graph of $du$ is invariant by the Hamiltonian flow.
\end{cor}
\demo
Let $u~: M\rightarrow \R$  satisfy the hypotheses of the corollary. We deduce from proposition \ref{propkam} that $u$ is a positive and negative weak KAM solution. Hence, $u$ is  $C^{1,1}$.  We have then~: 
$\forall t>0,\Gc (du)=  \overline{\Gc(dT_t^-u)}\subset \varphi_t(\Gc(du))$ and  $\Gc(du)= \overline{\Gc(dT_t^+u)}\subset \varphi_{-t}(\Gc(du))$). Hence the graph $\Gc(du)$ is invariant.
\enddemo

\section{Proof of   theorem \ref{th1}}\label{sectionproof}
 
Two submanifolds of $T^*M$ that are Hamiltonianly isotopic to the zero section have a non-empty intersection (see \cite{LaSi}). Let us now consider a submanifold $N$ of $T^*M$ that is Hamiltonianly isotopic to the zero section and that is invariant by the Tonelli flow of $H$. Then as $N$ is an invariant Lagrangian submanifold, there exists $k\in\R$ such that $N\subset \{ H=k\}$. Moreover, the intersection of $N$ with any $\Gc(du)$ for $u\in C^2(M, \R)$ is non-empty because the two manifolds are Hamiltonianly isotopic to the zero-section.  We have seen that Ma\~ n\'e's critical value is given by~: 
 $\displaystyle{c=\inf_{u\in C^\infty (M, \R)}\max_{q\in M}H(q, du(q))}$. Then we have~: $k\leq c$.
 
 We assume now that $N$ is a submanifold that is Hamiltonianly isotopic to the zero section and that is invariant under the Tonelli flow of $H$.  We have noticed that there exists $k\leq c$ such that $N\subset \{ H=k\}$. \\
Moreover, we have built in section \ref{select} a generating function $S$ and a function selector $\Phi$. There exists a dense open subset $U_0$ of $M$ with full Lebesgue measure such that $\Phi$ is differentiable on $U_0$ and~: $\forall q\in U_0, (q, d\Phi(q))\in N$. Hence, at Lebesgue almost every point, we have~: $H(q, d\Phi (q))\leq k$. We have seen that this implies~: $\Phi\prec L+k$.  As $k\leq c$ and $c$ is   the greatest lower bound of the set of the numbers $k\in \R$ for which there exists $u\in C^0(M, \R)$ with $u\prec L+k$, we deduce that $k=c$.
\subsection{Place of the Aubry set}
The beginning of this proposition is proved in \cite{Patpoltsib}.
\begin{prop} If $N$ is a submanifold that is Hamiltonianly isotopic to the zero section and that is invariant under the Tonelli flow of $H$, if $\Phi~: M\rightarrow \R$ is the associated function selector, then at every $q\in \Ac(H)$, $\Phi$ is differentiable, $(q, d\phi(q))\in N$ and $\Phi(q)=S\circ i_S^{-1}(q, d\Phi(q))$.
\end{prop}
We need a lemma~:
\begin{lemma}\label{nonsmooth}
Let $f~: U\rightarrow \R$ be a Lipschitz function defined on a open subset $U$ of $\R^d$ and let $U_0\subset U$ be a subset with full Lebesgue measure such that $f$ is differentiable at every point of $U_0$. We introduce a notation~: if $q\in U$, $K_f(q)$ is the  set of all the limits $\displaystyle{\lim_{n\rightarrow \infty} df(q_n)}$ where $q_n\in U_0$,  $\displaystyle{\lim_{n\rightarrow \infty}q_n=q}$ and $C_f(q)$ is the convex hull of $K_f(q)$. Then, at every point $q\in U$ where $f$ is differentiable, we have~: 
$df(q)\in C_f(q)$.
\end{lemma}

\demo This lemma is proved in \cite{FaMa}. A more general result is proved  in \cite{Cla1} too. Let us give an idea of a simple  proof. Using Fubini theorem, we obtain for every $v\in\R^d$ a sequence of vectors $(v_n)$ converging to $v$ and a decreasing sequence $(t_n)$ tending to $0$ such that~: 
$\displaystyle{df(q)v=\lim_{n\rightarrow \infty}\frac{1}{t_n}\int_0^{t_n}df(q+sv_n).v_nds}$ where for  Lebesgue almost every point $t\in [0, t_n]$, we have~: $q+tv_n\in U_0$. Then, for every $v\in \R^d$, we find $p_v\in C_f(q)$ such that $df(q)v=p_v(q)$; as $C_f(q)$ is convex and compact, using Hahn-Banach theorem, we deduce~: $df(q)\in C_f(q)$.
\enddemo
As $\Phi\prec L+c$, there exists a pair $(u_-, u_+)$ of conjugate weak KAM solutions such that $u_{-|\Mc(H)}=u_{+|\Mc(H)}=\Phi_{|\Mc(H)}$ and we have~: $u_+\leq \Phi\leq u_-$. As $u_+$ and $u_-$ are differentiable on $\Ac(H)$ and as $u_{-|\Ac(H)}=u_{+|\Ac(H)}=\Phi_{|\Ac(H)}$ and $du_{-|\Ac(H)}=du_{+|\Ac(H)}$, we deduce that $\Phi$ is differentiable on $\Ac(H)$ and that~: $\forall q\in \Ac(H), (q, d\Phi (q))=(q, du_-(q))\in \Ac^*(H)$. We cannot conclude that $\Ac^*(H)\subset N$ because we don't know if $\Ac(H)\subset U_0$. We use then lemma \ref{nonsmooth} (we work in a chart). Let $q_0\in \Ac(H)$ be an element of the projected Aubry set. We deduce from the lemma that~: $d\Phi(q_0)\in C_\Phi(q_0)$. Moreover, $d\Phi(q_0)\in T_{q_0}^*M\cap \{ H=c\}$ and $T_{q_0}^*M\cap \{ H=c\}$ is the set of the extremal points of the convex set  $T_{q_0}^*M\cap \{ H\leq c\}$ and this last set contains $C_\Phi(q_0)$. Then $d\Phi(q_0)$ is an extremal point of $C_\Phi(q_0)$ and then $d\Phi(q_0)$ belongs to $K_\Phi(q_0)$. It means that there exists a sequence $(q_n)$ of points of $U_0$ that converge to $q_0$ so that~: $\displaystyle{d\Phi(q_0)=\lim_{n\rightarrow \infty} d\Phi(q_n)}$. We deduce that $(q_0, d\Phi(q_0))\in N$. Moreover, as $\Phi$, $S$ and $i_S$  are continuous~: $$\Phi(q_0)=\lim_{n\rightarrow \infty} \Phi(q_n)=\lim_{n\rightarrow \infty} S\circ i_S^{-1}(q_n, d\Phi(q_n))=S\circ i_S^{-1}(q_0, d\Phi(q_0)).$$
 \subsection{ Place of the non-wandering set}
\begin{prop}  If $N$ is a submanifold that is Hamiltonianly isotopic to the zero section and that is invariant under the Tonelli flow of $H$, we have~: $\Omega(\varphi_{t|N})\subset \Ac^*(H)$.

\end{prop}
 
Let us  explain why the non-wandering set $\Omega (\varphi_{t|N})$ of the Hamiltonian flow restricted to $N$ is in the Aubry set for $H$. A similar argument is given is \cite{Patpoltsib}.   We have noticed that $N\subset \{ H=c\}$. Let $(q, p)\in \Omega (\varphi_{t|N})$. Then there exist a sequence $(q_n, p_n)$ of points of $N$ converging to $(q,p)$ and a sequence $(t_n)$ in $\R_+$ tending to $+\infty$ such that~: $\displaystyle{\lim_{n\rightarrow \infty} \varphi_{t_n}(q_n, p_n)=(q,p)}$. Let us introduce the notation~: $(q_n(t), p_n(t))=\varphi_t(q_n, p_n)$.   As $N$ is Hamiltonianly isotopic to the zero-section, it is exact Lagrangian and then~: $\displaystyle{\lim_{n\rightarrow \infty} \int_0^{t_n}p_n(t)\dot q_n(t)dt=0}$. As $(q_n(t), p_n(t))$ is an orbit, we have~: $\dot q_n=\frac{\partial H}{\partial p}(q_n, p_n)$, and then~:  
 $p_n.\dot q_n=L(q_n, \dot q_n)+H(q_n, p_n)= L(q_n, \dot q_n)+c$.  Finally, we have~:
 $$\lim_{n\rightarrow \infty} \int_0^{t_n}(L(q_n(t), \dot q_n(t))+c)dt=0.$$
 As $\displaystyle{\lim_{n\rightarrow \infty}q_n(0)=q}$ and  $\displaystyle{\lim_{n\rightarrow \infty}q_n(t_n)=q}$, we deduce that~: 
$$(q, p)=\lim_{n\rightarrow \infty} (q_n(0), p_n(0))=\lim_{n\rightarrow \infty}(q_n(0), \frac{\partial L}{\partial v}(q_n(0), \dot q_n(0)))\in \Ac^*(H).$$
  Hence we have proved~:
 $$ \Omega (\varphi_{t|N})\subset \Ac^* (H).$$
 \subsection{Comparison between $h$ and $\Phi$}
  Let us now prove that $\Phi$ satisfies the hypotheses of corollary \ref{corkam}.
 \begin{prop}
 For all $q\in U_0$, there exists $q_1, q_2\in \Ac(H)$  such that~: $\Phi(q_1)-\Phi(q)\geq h(q, q_1)$ and $\Phi(q)-\Phi(q_2)\geq h(q_2, q)$.
 \end{prop}
 
 We consider $q\in U_0$.  Then  $(q, d\Phi (q))\in N$ and $\Phi (q)=S\circ i_S^{-1}(q, d\Phi(q))$.  Then the $\alpha$ and $\omega$ limit sets of $(q, d\Phi(q))$ are non-empty. There exist $(q_1, p_1)\in \omega (q, d\Phi(q))$ and $(q_2, p_2)\in \alpha(q, d\Phi(q))$. These points being non-wandering and in $N$, we have noticed that they belong to $\Ac^*(H)$~: $(q_i, p_i)=(q_i, d\Phi(q_i))\in \Ac^*(H)$ and that~: $\Phi(q_i)= S\circ i_S^{-1}(q_i, d\Phi(q_i))$.  As they belong to the $\alpha$/$\omega$ limit set, there exist two sequences $(t_n)\in\R_+$ and $(\tau_n)\in\R_+$  tending to $+\infty$ so that~:
 $$\lim_{n\rightarrow \infty}\varphi_{t_n}(q, d\Phi(q))=(q_1, p_1) \quad{\rm and}\quad \lim_{n\rightarrow \infty}\varphi_{-\tau_n}(q, d\Phi(q))=(q_2, p_2).$$
 We use the following notation~: $\varphi_t(q, d\Phi(q))=(q(t), p(t))$ and we compute~:
 
$$\Phi(q_1)-\Phi(q)=S\circ i_S^{-1}(q_1, d\Phi (q_1))-S\circ i_S^{-1}(q, d\Phi(q))=\lim_{n\rightarrow \infty} S\circ i_S^{-1}\circ \varphi_{t_n}(q, d\Phi(q))-S\circ i_S^{-1}(q, d\Phi(q))$$
We have~: $$S\circ i_S^{-1}\circ \varphi_{t_n}(q, d\Phi(q))-S\circ i_S^{-1}(q, d\Phi(q))=S\circ i_S^{-1}(q(t_n), p(t_n))-S\circ i_S^{-1}(q(0), p(0))$$
where~: $i_s(q, \zeta)=(q, \frac{\partial S}{\partial q}(q, \zeta))$ and on $\Sigma_S$~: $\frac{\partial S}{\partial \zeta}=0$.  Hence $i_s^{-1}(q,p)=(q, \beta (q,p))$ and~: $\forall ( \delta q, \delta p)\in T_{(q,p)}N$~: $d(S\circ i_S^{-1})(q,p)(\delta q, \delta p)=\frac{\partial S}{\partial q}(i_S^{-1}(q, p))\delta q$. Then~: 
$$S\circ i_S^{-1}\circ \varphi_{t_n}(q, d\Phi(q))-S\circ i_S^{-1}(q, d\Phi(q))=\int_0^{t_n}\frac{\partial S}{\partial q}(i_S^{-1}(q(t), p(t)))\dot q(t)dt=\int_0^{t_n} p(t)\dot q(t)dt$$
As $(q(t), p(t))$ is an orbit, we have~: $\dot q=\frac{\partial H}{\partial p}(q, p)$, and then~:  
 $p.\dot q=L(q, \dot q)+H(q, p)= L(q, \dot q)+c$.  Finally, we have~:

$$S\circ i_S^{-1}\circ \varphi_{t_n}(q, d\Phi(q))-S\circ i_S^{-1}(q, d\Phi(q)) =\int_0^{t_n}(L(q(t), \dot q(t))+c)dt.$$
We deduce that~:
$$h(q, q_1)\leq \lim_{n\rightarrow \infty} S\circ i_S^{-1}\circ \varphi_{t_n}(q, d\Phi(q))-S\circ i_S^{-1}(q, d\Phi(q))=\Phi(q_1)-\Phi(q).$$
In a similar way, we obtain~: $h(q_2, q)\leq \Phi(q)-\Phi(q_2)$. 
\subsection{Conclusion}

We deduce from this and from corollary \ref{corkam} that $\Phi$ is $C^{1,1}$ and that $\Gc(d\Phi)$ is invariant by the flow. 

Let us now summarize what we did~:\begin{enumerate}
\item[$\bullet$] we have found a dense part $\Gc(d\Phi_{|U_0})$ of $\Gc(d\Phi)$ that is a subset of the closed manifold $N$. Hence $\Gc (d\Phi)\subset N$;
\item[$\bullet$] hence $\Gc (d\Phi)$ is a closed submanifold of $N$ that has the same dimension as $N$; $N$ being connected, we deduce that  $\Gc(d\Phi)=N$ is a graph.
\end{enumerate}

\end{document}